\documentclass[11pt]{article}
\usepackage{graphicx}
\usepackage{amsmath,amssymb,amsthm,amsfonts}
\usepackage[normal]{subfigure}
\usepackage{amssymb}
\usepackage{subfig}
\usepackage{amscd}
\usepackage{verbatim}
\newcommand{\vseven}{\\[7pt]}
\newcommand{\hf}{\hfill}
\newcommand{\hsf}{\hspace*{\fill}}

\newtheorem{theorem}{Theorem}[section]

\newtheorem{lemma}[theorem]{Lemma}
\newtheorem{proposition}[theorem]{Proposition}
\newtheorem{remark}[theorem]{Remark}

\newcommand{\R}{\mbox{$\mathbb{R}$}}

\newcommand{\Rn}{\mbox{$\mathbb{R}^N$}}
\newcommand{\RN}{\mbox{$\mathbb{R}^N$}}
\newcommand{\N}{\mbox{$\mathbb{N}$}}
\newcommand{\Z}{\mbox{$\mathbb{Z}$}}

\newcommand{\ag}{\alpha}

\newcommand{\dg}{\delta}

\newenvironment{acknowledgement}{\textbf{Acknowledgement.}\em}{}
\def\endproof{\hfill$\square$\vspace{6pt}}
\begin{document}
\title{On a class of singular second-order Hamiltonian systems with infinitely many homoclinic solutions}
\author{David G. Costa, Hossein Tehrani\\
\begin{tabular}{l}
{\small Dept of Math Sciences, University of Nevada Las Vegas, Box 454020, Las Vegas, NV 89154-4020, USA}
\end{tabular}\\
\begin{tabular}{c}
{\small costa@unlv.nevada.edu, tehranih@unlv.nevada.edu}
\end{tabular}}
\date{}

\maketitle
{\small \textbf{Abstract:}
We show existence of infinitely many homoclinic orbits at the origin for a class of singular second-order Hamiltonian systems
$$
\ddot{u} + V_u (t,u)=0\,,\quad -\infty < t < \infty\,.
$$
We use variational methods under the assumption that\ $V(t,u)$\ satisfies the so-called "Strong-Force" condition.\\[0.15cm]

{\small \textbf{Resum\'{e}:}
Nous prouvons l'existence d'un nombre infini d'orbites homocliniques \`{a} l'origine pour une classe de  syst\`{e}mes Hamiltoniens singuliers du second ordre,
$$
\ddot{u} + V_u (t,u)=0\,,\quad -\infty < t < \infty\,.
$$
Nous utilisons les m\'{e}thodes variationelles avec l'hypoth\`{e}se que $V(t,u)$ satisfait la condition dite de "Strong-Force".

\section*{Introduction}                                             
The search for periodic as well as homoclinic and heteroclinic solutions of Hamiltonian systems has a long and rich history. In this paper we are particularly interested in homoclinic solutions of singular second-order Hamiltonian systems with time periodic potentials. We refer the interested reader to the book \cite{1} of Ambrosetti and Coti Zelati for results on the literature of periodic solutions for such singular systems.

Second-order Hamiltonian systems are systems of the form
\vseven
$(HS)$\hf $\displaystyle \ddot{u} + V_u (t,u)=0\,\quad\mbox{ $t\in\R$, $u\in\RN$}\, .$\hsf
\vseven

Loosely speaking, they are the Euler-Lagrange equations of the functional
$$
I(u)=\int_{}L(t,u,\dot{u})\,dt\,,
$$
where the integration is taken over a finite interval $[0,T]$ or all reals $\R$ and the {\it Lagrangian} has the form
$$
L(t,u,\dot{q})=\frac{1}{2}|\dot{u}|^2 - V(t,u)\,.
$$
Clearly, when the potential $V(t,u)$ is $T$-periodic in $t$, it is natural to look for $T$-periodic solutions of $(HS)$ as critical points of the functional $I(u)$ over a suitable space of $T$-periodic functions. Also, in such a case, one can look for homoclinic solutions at the origin (i.e., solutions of $(HS)$ satisfying $u(t)$, $\dot{u}(t)\longrightarrow 0$) as limits of $kT$-periodic solutions ({\it subharmonic solutions}) as $k\to\infty$ (see \cite{16}) or, alternatively, as critical points of the functional $I(u)$ over a suitable space of functions on the whole space $\R$ (typically, $H^1 (\R,\RN)$).

For singular systems, one assumes that $V\in C^1(\R\times \R^N\setminus S)$ and $\lim_{u\rightarrow S}|V(t,u)|=\infty$ for some $S\subset \R^N$.  Although the study of singular systems is perhaps as old as the Kepler's classical problem in mechanics,
$$
\ddot{u} + \frac{u}{|u|^3}=0
$$
(and, also, the $N$-body problem), the interest in such problems was renewed by the pioneering papers \cite{12} of Gordon in 1975 and \cite{13} of Rabinowitz in 1978. In \cite{12} the notion of Strong-Force is introduced to deal with singular problems, while in \cite{13} the use of variational methods is brought into the study of periodic solutions of Hamiltonian systems.

The present paper is concerned with existence of homoclinic solutions for second-order Hamiltonian systems
\vseven
$(HS)$\hf $\displaystyle \ddot{u} + V_u (t,u)=0\,,$\hsf
\vseven
where\ $-\infty < t < \infty$, $u=(u_1, u_2,\ldots ,u_N)\in\Rn$\ and the potential\ $V:\R\times \R^N\setminus \{q\}\rightarrow \R$\ has a singularity $0\neq q\in\Rn$. We recall that a homoclinic solution of\ $(HS)$\ is a solution such that $u(t)\in R^N\setminus \{q\}$ for all $t\in \R$ and
$$
u(t),\dot{u}(t)\longrightarrow 0\quad\mbox{as $t\to\pm\infty$\,.}
$$
Throughout the paper we will be considering the following assumptions on\ $V(t,u)$:
\begin{itemize}
\item[$(A)$]\ $V(t,u)=a(t)W(u)$, with $a\in C(\R)$ a $T$-periodic function such that $a_0 \leq a(t)\leq a_{\infty}$\ \mbox{ for some\ $a_0, a_{\infty} > 0$}\,;
\item[$(H_1)$]\ $\displaystyle$ $W\in C^2 (\Rn\setminus\{ q\},\R)$\ \mbox{ for some\ $q\in\Rn\setminus\{ 0\}$}\,;
\item[$(H_2)$]\ $\displaystyle$ \mbox{$W(0)=W_u (0)=0$},\ $W(u) < W(0)=0$\ for\ $u\neq 0$,\ and\ $-\ag_0 I \leq W_{uu}(0)\leq -\ag_1 I$\ \mbox{ for some\ $\ag_0, \ag_1 > 0$}\,;
\item[$(H_3)$]\ $\displaystyle$ $\lim_{u\to q}W(u)= -\infty$\ and there exists\ $U\in C^1 (\Rn\setminus\{ q\},\R)$ such that\ $\lim_{u\to q}|U(u)|=\infty$\ and\ $W(u)\leq -|\nabla U(u)|^2$\ for\ $0<|u-q|\leq r$\,;
\item[$(H_4)$]\ $\displaystyle$ There exists\ $U_{\infty}\in C (\Rn\setminus B_{R_0},\R)$\ such that\ $\lim_{|u|\to\infty}|U_{\infty}(u)|=\infty$\ and\ $W(u)\leq -|\nabla U_{\infty}(u)|^2$\ for\ $u$\ large\,.
\end{itemize}
\medskip

Note that by our assumptions, $W$ has a strict global maximum at $u=0$ which by $(H_2)$ is an unstable equilibrium of\ $(HS)$.  Furthermore \ $(H_3),(H_4)$\ concern the behavior of $W$ close to the singularity and at infinity. In fact, $(H_3)$ indicates that the potential $W$ satisfies the {\it Strong-Force} condition mentioned earlier (used by Gordon in \cite{12}) which governs the rate at which $W(x)$ approaches $-\infty$ as $x\rightarrow q$. A typical example is\ $W(x) = |x - q|^{-\alpha}$\ ($\alpha\geq 2$)\ in a neighborhood of\ $q$. On the other hand $(H_4)$ allows $W$ to go to zero at infinity although at a slow rate. This condition will be satisfied if, for example, $\lim_{|x|\rightarrow\infty}|x|^{\beta}W(x)\neq 0$ for some $\beta\in (0,2]$.

In the case of autonomous singular Hamiltonian systems, the first result on existence of a homoclinic orbit using variational methods was obtained by Tanaka \cite{20} under essentially the same assumptions as above.  In \cite{20} Tanaka used a minimax argument from\ Bahri-Rabinowitz \cite{2}\ in order to get approximating solutions of the boundary value problems
$$
\ddot{u} + V^{\prime} (u)=0\,,\quad t\in (-m,m)\,,\quad u(-m)=u(m)=0
$$
as critical points of the corresponding functionals, and obtained uniform estimates to show that those solutions converged weakly to a {\it nontrivial} homoclinic solution of\ $(HS)$.  Later Bessi \cite{4},
using Lusternick-Schnirelman category, proved the existence of $N-1$ distinct homoclinics for potentials satisfying a pinching condition (see also \cite{1} and \cite{21} for multiplicity results in case of smooth Hamiltonians). Different kinds of multiplicity results were obtained in \cite{3,6} (still for conservative systems) by exploiting the topology of $\R^N\setminus S$, the domain of the potential, when the set $S$ is such that the fundamental group of $\R^N\setminus S$ is nontrivial.

In the case of planar autonomous systems more extensive existence and multiplicity results were obtained.
Indeed, under essentially the same conditions as above with\ $N=2$, Rabinowitz showed in \cite{15} that $(HS)$ has at least a pair of homoclinic solutions by exploiting the topology of the plane and minimizing the energy functional on classes of sets with a fixed winding number around the singularity $q$ (see also \cite{5} for results in the case of two singularities). The result in \cite{15} was substantially improved in \cite{7} where, using the same idea, the authors show that a nondegeneracy variational condition introduced in \cite{15} is in fact necessary and sufficient for the minimum problem to have a solution in the class of sets with winding number greater than 1 and, therefore, proved  a result on existence of infinitely many homoclinic solutions.

On the other hand, in the case of $T$-periodic time dependent Hamiltonians in $\R^N$, existence of infinitely many homoclinic orbits was obtained for {\it smooth} Hamiltonians by using a variational procedure due to S\'{e}r\'{e} in \cite{17} and \cite{18} for first order systems, and in \cite{8} and \cite{9} for second order equations. In the case $N=2$, using these ideas, Rabinowitz \cite{14} constructed infinitely many multibump homoclinic solutions for\ $V(t,u)$\ of the form\ $a(t) W(u)$, with\ $a(t)$\ being almost periodic and\ $W(u)$\ satisfying\ $(A),(H_1)$-$(H_4)$.

Our work here on homoclinic solutions of time periodic singular equations was motivated by earlier works on periodic solutions of such equations as well as by \cite{9}, where homoclinic solutions in $\R^N$ are considered in the case of second order {\it smooth} Hamiltonians. As was mentioned above, the main tool in \cite{9} is a minimax procedure of S\'{e}r\'{e} which gives the existence of infinitely many multibump homoclinics.
The novelty in our approach is the use of category theory in the case of homoclinics. Let us now indicate the main steps in our approach. In order to find homoclinics, we consider the action functional $I$ on the full space $\Lambda=H^1(\R, \R^N\setminus \{q\})$ but, as in the periodic case, use  Lusternick-Schnirelman category theory to generate a sequence of minimax values that are candidates for critical levels of $I$. Indeed, since $Cat_{\Lambda}(\Lambda)=\infty$\,(cf. Proposition \ref{Cat}), this process can be initiated and we can define the sequence
$$
c_k:=\inf_{S\in\Gamma_k}\sup_{u\in S}I(u)\,\quad k\in\mathbb{N},
$$
where
$$
\Gamma_k := \{\, S\subset\Lambda\, |\, S \mbox{ is compact and }Cat_{\Lambda}(S)\geq k\,\}\,.
$$
By contrast with the case of periodic solutions, the homoclinic problem exhibits a lack of compactness (indeed Palais-Smale condition is not satisfied) which makes the application of critical point theorems quite challenging. However, we will show that there is enough control on the Palais-Smale sequences to prove the existence of one homoclinic solution as the weak limit of a (PS)-sequence corresponding to the first positive minimax level\, $c_2$ above. This is done in the section 1. We point out that the approach in section 1 is applicable to autonomous systems as well and, therefore, provides an independent proof of the result of Tanaka mentioned above. In the second section, by adapting some of the ideas in \cite{9}, we will show how to get a complete description of the behavior of (PS)-sequences and a suitable version of a deformation theorem, thus allowing application of variational methods. Finally, by using these results and an indirect argument, we prove that, under conditions $(A)$, $(H_1)-(H_4)$, the singular second-order Hamiltonian system
$$
\ddot{u} + a(t) W^{\prime}(u)=0\,,\quad -\infty < t < \infty\,,
$$
possesses infinitely many {\it geometrically distinct} homoclinic solutions.

To the best of our knowledge, this is the first result on existence of infinitely many homoclinics for time periodic and singular second-order Hamiltonian systems in $\R^N$\ ($N\geq 3$) when the singularity is a point.
\section{Existence of a nontrivial homoclinic solution}                  
The main result of this section is the following:
\setcounter{theorem}{0}
\begin{theorem}\label{theorem1.1}
Assume conditions $(A)$, $(H_1)$-$(H_4)$ stated in the Introduction. Then, the singular Hamiltonian system
\vseven
$(SHS)$ \hf
$ \displaystyle \ddot{u} + a(t) W^{\prime}(u)=0\,,\quad -\infty < t < \infty\,,$  \hsf
\vseven
has at least one homoclinic solution emanating from zero.
\end{theorem}

The rest of this section is devoted to a proof of this result.
Let\ $H^1=H^1 (\R,\Rn)$\ denote the usual Sobolev space with inner-product
$$
\langle u , w\rangle = \int_{-\infty}^{\infty}[\dot{u}\cdot\dot{w} + u\cdot w]\,dt\quad\ \forall u,w\in H^1
$$
and corresponding norm\ $||u||=\langle u , u\rangle^{\frac{1}{2}}$. We consider the open set
$$
\Lambda := \{\,u\in H^1\ |\ u(t)\neq q\ \forall t\in\R\,\}
$$
and the functional given by
\begin{equation}\label{functional}
I (u) := \frac{1}{2}\int_{-\infty}^{\infty}|\dot{u}|^2\,dt - \int_{-\infty}^{\infty}a(t)W(u)\,dt\,,\quad u\in\Lambda\,.
\end{equation}
Critical points of $I$ are solutions of the Hamiltonian system $(SHS)$ and, as we shall see below, our hypotheses will guarantee that any such critical point is in fact a homoclinic solution emanating from zero.  We start by stating, in the next two lemmas, some basic properties of the space $H^1$ and the sublevel sets of the functional $I$ which will be used throughout the presentation. But first a word on notation: unless otherwise indicated, we assume that all integrals are taken over $\R$ . Furthermore, for $1\leq p\leq\infty$, we denote the norm of the corresponding $L^p(\R,\R^N)$  space by $||\cdot||_p$.
\begin{lemma}\label{lemma1.0}
If\, $v\in H^1=H^1 (\R,\Rn)$\, then $v\in C^{0,\frac{1}{2}}(\R)$ and, for $s\in \R$,
\begin{equation}\label{eq1}
|v(s)| \leq \left( \int_{A(s)}|v|^2(t) dt \right)^{\frac{1}{2}}+\left( \int_{A(s)}|\dot{v}|^2(t) dt \right)^{\frac{1}{2}}\,,
\end{equation}
where
$$A(s)=\left\{
\begin{array}{cc}
[s,s+1] & \mbox{ when } s\geq 0  \\

[s-1,s] & \mbox{ when } s <0\,.
\end{array}
\right.$$
\end{lemma}
\begin{lemma}\label{lemma1.0.0}
(cf. Lemma 2.2 in \cite{7}) Consider the above functional $I$ defined on the whole space $H^1$, so that
$I:H^1\rightarrow [0,\infty]$. Then $I$ is weakly lower semicontinuous. Furthermore, for $b\in\R^{+}$, denote $I^{b}=\{ u\in H^1: I(u)\leq b\}$. Then
\begin{description}
\item[(a)] There exists $R=R(b) >0$ such that, $||u||_{\infty}\leq R$\  for all $u\in I^{b}$.
\item[(b)] There exists $\rho=\rho(b) >0$ such that $dist(range (u), q)\geq \rho$\ for all   $u\in I^{b}$.
\item[(c)] For given $\dg > 0$, there exists $\tau_{\delta}=\tau_{\delta}(b)>0$ such that $meas(S_{\delta}(u))\leq\tau_{\delta}$\,,where\ $S_{\delta}(u):=\{\, t\in\R\,|\, |u(t)|\geq\delta\,\}$\,.
\end{description}
\end{lemma}
We present below further properties of the functional $I$ on $\Lambda$ which will enable us to set up a variational characterization of some of its critical values.
\begin{lemma}\label{lemma1.1}
$I:\Lambda\longrightarrow\R$ is well-defined and of class $C^1$.
\end{lemma}
\proof  In view of Lemma \ref{lemma1.0}, any $u\in\Lambda\subset H^1$ is H\"older continuous of exponent $\frac{1}{2}$ and $\lim_{|t|\to\infty}u(t)=0$, so that
\begin{itemize}
\item[$(i)$]\ Given $\epsilon > 0$, there exists $T_{u,\epsilon}>0$ such that $|u(t)|\leq\epsilon$ for $|t|\geq T_{u,\epsilon}$\,;
\item[$(ii)$]\ There exists $r=r(u)>0$ such that $|u(t) - q|\geq r$\,\ $\forall t\in\R$.
\end{itemize}
Now, by $(H_2)$, there exists $\dg > 0$ such that
\begin{equation}\label{eq1}
\frac{\ag_0}{2}|u|^2 \leq - W(u) \leq 2\ag_1 |u|^2\quad\mbox{ if }\ |u|\leq\dg\,.
\end{equation}
Therefore, for a given $u\in \Lambda$ (with $T_{\delta}=T_{u,\delta}>0$ and $r=r(u)$ as in $(i),(ii)$ above), we have
$$
-\int_{}a(t)W(u)\,dt \leq -\int_{|t|\leq T_{\delta}}a(t)W(u)\,dt -\int_{|t|\geq T_{\delta}}a(t)W(u)\,dt \leq 2a_{\infty}T_{\delta}M_1 +  2a_{\infty}\alpha_1\int_{}|u|^2dt
$$
where $M_1 =\max_{}\{\,|W(z)|\ |\ |z|\leq\| u\|_{\infty},\ dist (z,u)\geq r\, \}$, so that
$$
-\int_{}a(t)W(u)\,dt \leq C_1\| u\|^2_2 + C_2\,.
$$
This shows that the functional $I$ is well-defined on $\Lambda$. In a similar manner, by using $(i),(ii)$ and $(H_2)$, one can show that $I$ is of class $C^1$ on $\Lambda$\, and
$$
\left< I'(u),\phi \right> =\int \dot{u}\dot{\phi}dt -\int a(t) W'(u)\phi dt\quad\quad\forall u, \phi\in H^1\,.
$$
\endproof
\begin{lemma}\label{lemma1.2}
$I:\Lambda\longrightarrow\R$ is coercive, i.e., if $u_n\in\Lambda$ is such that $I(u_n)\leq b$ for some $b>0$ then $\| u_n\|$ is bounded by a constant depending only on $b$.
\end{lemma}
\proof
Let $u_n\in\Lambda$ be such that $I(u_n)\leq b$ for some $b>0$. Taking $\delta > 0$ as in \eqref{eq1}, we get
$$
-\int_{}a(t)W(u_n)\,dt = -\int_{|u_n (t)|\leq \delta}a(t)W(u_n)\,dt -\int_{|u_n (t)|\geq \delta}a(t)W(u_n)\,dt \geq a_0\int_{|u_n (t)|\leq \delta}\frac{\alpha_0}{2}|u_n|^2\,dt
$$
Furthermore, using Lemma \ref{lemma1.0.0},
$$
\int_{|u_n (t)|\geq \delta}|u_n|^2\,dt \leq R(b)^2 meas(S_{\delta}(u_n)) \leq R(b)^2 \tau_{\delta}(b) := C_{\delta}(b)\,,
$$
where $C_{\delta}=C_{\delta}(b)>0$ is independent of $n$\,. The above two estimates give
$$
\frac{a_0 \ag_0}{2}\int_{}|u_n|^2\,dt \leq - \int_{}a(t)W(u_n)\,dt + \frac{a_0 \ag_0}{2} C_{\delta}(b)  \leq I(u_n) + \hat{C}_{\delta}\,.
$$
Since boundedness of $I(u_n)$ implies that $\| \dot{u}_n\|_2^2$ is bounded, we conclude that $\| u_n\|^2 = \| \dot{u}_n\|_2^2 + \| u_n\|_2^2$ is bounded as well.
\endproof

We recall that a sequence $(u_n)\subset \Lambda$ is called a $(PS)_c$-sequence for $I$ if $I(u_n)\rightarrow c$ and $ I'(u_n)\rightarrow 0 $ as $n\rightarrow\infty$. In addition, $I$ is said to satisfy the $(PS)_c$ condition if any $(PS)_c$-sequence has a convergent subsequence (to a critical point of $I$). As is well-known, some version of such a compactness assumption is necessary for application of critical point theorems.  A complete description of the behavior of $(PS)_c$-sequences of $I$  will be given in Theorem 2.1 of the next section. The following (weaker) result is all that we need here.

\begin{lemma}\label{lemma1.3}
If $(u_n)$ is a $(PS)_c$-sequence for $I$ for some $c>0$, then there exists another $(PS)_c$-sequence $(v_n)$ such that $v_n \rightharpoonup v$ in $H^1$ for some nonzero $v\in\Lambda$.
\end{lemma}
\proof
In view of Lemma \ref{lemma1.2}, there exists $u\in H^1$ such that $u_n \rightharpoonup u$ and $u_n (t)\rightarrow u(t)$ {\it locally uniformly} on $\R$ (i.e., uniformly on compact subsets of $\R$). Since $I$ is weakly lower-semicontinuous, it follows that $I(u)\leq c$ and, hence, $u\in\Lambda$ (see Lemma \ref{lemma1.0.0}).

Now, let $t_n \in\R$ be such $|u(t_n)|=\max_{t\in\mathbb{R}}|u(t)|$ and define
$$
v_n (t) = u_n (t + l_n T)\,,
$$
where $l_n\in\Z$ satisfies $l_n T \leq t_n < (l_n + 1)T$. Then, there exists $\hat{t}_n \in [0,T)$ such that
$$
v_n (\hat{t}_n)=\max_{t\in\mathbb{R}}|v_n (t)|=\max_{t\in\mathbb{R}}|u_n (t)|\,,
$$
and, from the $T$-periodicity of $a(t)$, we can easily check that
$$
I(v_n) = I(u_n)\ \mbox{ and }\ \| I^{\prime}(v_n)\| = \| I^{\prime}(u_n)\|\,.
$$
Therefore, $(v_n)$ is also a $(PS)_c$-sequence and we may assume that
$$
v_n \rightharpoonup v\ \mbox{ and }\ v_n (t)\rightarrow v(t)\,\ \mbox{{\it locally uniformly on $\R$}}\,.
$$
It remains to show that $v\neq 0$ (of course  $v\in\Lambda$).
\medskip

Indeed, if $v\equiv 0$, then $v_n \to 0$ uniformly on $[0,T]$ and, consequently, $v_n \to 0$ uniformly on $\R$ since $\max_{t\in\mathbb{R}}|v_n (t)| = |v_n (\hat{t}_n)|\longrightarrow 0$ as $n\to\infty$.

On the other hand, $(H_2)$ implies that
\begin{equation}\label{eq3}
-W(u)\geq \frac{\ag_0}{2}|u|^2\,\ \mbox{ and }\ - W^{\prime}(u)\cdot u = - W^{\prime\prime}(0)u\cdot u + o(|u|^2) \geq \frac{\ag_0}{2}|u|^2
\end{equation}
for $|u|$ small. Moreover, since $(v_n)$ is a bounded $(PS)$-sequence, we have
\begin{equation}\label{eq4}
o(1)=o(1)\| v_n\| = I^{\prime}(v_n).v_n = \| \dot{v}_n\|_2^2 - \int_{}a(t)W^{\prime}(v_n)v_n\,dt
\end{equation}
where\ $\| \dot{v}_n\|_2^2$\ is bounded and, in view of \eqref{eq3} and $(A)$,
\begin{equation}\label{eq5}
\left| -\int_{}a(t)W^{\prime}(v_n)v_n\,dt \right| \geq \frac{a_0\ag_0}{2}\|v_n \|_2^2\,.
\end{equation}
Therefore, \eqref{eq4} and \eqref{eq5} give
$$
\| v_n\|^2 = \| v_n\|_2^2 + \|\dot{v}_n\|_2^2 = o(1)\,,
$$
which contradicts the fact that $I(v_n)\longrightarrow c > 0$\,.
\endproof

Next we plan to use Lusternik-Schnirelman (LS) category theory to construct a critical level for $I$. But first we will present a topological result that makes such a construction possible. In what follows, $Cat_{\Lambda}(Y)$ denotes the (LS) category of $Y\subset \Lambda$ with respect to $\Lambda$.
\begin{proposition}\label{Cat}
The following results hold:
\begin{enumerate}
\item  $Cat_{\Lambda}(\Lambda)=\infty$.
\item  For every $k\geq 1$, there exists a compact $Y\subset \Lambda$, such that $Cat_{\Lambda}(Y)=k$.
\end{enumerate}
\end{proposition}

We defer the proof of this proposition to the Appendix.  Using proposition \ref{Cat}, we define a sequence of minimax
values as follows:
$$
c_k:=\inf_{S\in\Gamma_k}\sup_{u\in S}I(u)\,,
$$
where
$$
\Gamma_k := \{\, S\subset\Lambda\, |\, S \mbox{ is compact and }Cat_{\Lambda}(S)\geq k\,\}\,.
$$
Note that $c_1=0$ (since $Cat_{\Lambda}(Y)=1$ for any singleton $Y\subset \Lambda$). Our next task is to show that $c_2>0$, and then provide a nonzero critical point of $I$.

\begin{lemma}\label{lemma1.4}
For any $S\in \Gamma_2$ (i.e. $S\subset \Lambda$ compact with $Cat_{\Lambda}(S)\geq 2$) there exist $v\in S$, $t_0 \in\R$, $k_0 > 1$ such that
\begin{equation}\label{eq6}
v(t_0)=k_0 q\,.
\end{equation}
\end{lemma}
\proof
Assume that \eqref{eq6} does not hold for any $u\in S$, i.e., no $u\in S$ has a range intersecting the half-line $\{\, tq\, |\, t>1\, \}$. Then we can define the linear homotopy $H:[0,1]\times S\longrightarrow \Lambda$ given by
$$
H(\tau,u) = (1-\tau) u\,,
$$
which shows that $Cat_{\Lambda}(S)=1$, a contradiction.
\endproof

Now, consider the subset $E\subset\Lambda$ defined by
$$
E := \{\, v\in\Lambda\, |\, \exists t_0\in\R, k_0 > 1 \mbox{ such that } v(t_0)=k_0 q\, \}\,.
$$
Note that, by Proposition \ref{Cat} and the above lemma, $E$ is nonempty.
\medskip

\begin{lemma}\label{lemma1.5}
The infimum $d:=\inf_{v\in E}I(v)$ is achieved (so, in particular, $d>0$)
\end{lemma}
\proof
Let $v_n \in E$ be a minimizing sequence for $d$, i.e.,
$$
I(v_n)\longrightarrow d\ \mbox{ as }\ n\to\infty
$$
and there exist $t_n\in\R$, $k_n > 1$ such that
$$
v_n (t_n) = k_n q\,.
$$
Since $\| v_n\|$ is bounded by Lemma \ref{lemma1.2}, it follows that $\| v_n\|_{\infty}$ and $k_n > 1$ are also bounded. Now, as in the proof of Lemma \ref{lemma1.3}, letting $l_n\in\Z$ satisfy $l_n T \leq t_n < (l_n + 1)T$ and defining $w_n (t) := v_n (t + l_n T)$, we have that
$$
I(w_n) = I(v_n)\longrightarrow d
$$
and
$$
w_n (s_n) = k_n q\,,
$$
where $s_n := t_n - l_n T \in [0,T)$. Thus, $(w_n )$ is another minimizing sequence for $d$ in $E$. And since $\|w_n\|_{\infty}=\| v_n\|_{\infty}$ and $k_n > 1$ are bounded, we have (passing to a subsequence, if necessary) that\ $w_n \rightharpoonup w$, $w_n (t) \longrightarrow w(t)$ {\it locally uniformly} on $\R$, $s_n \longrightarrow \hat{s}\in [0,T]$ and $k_n \longrightarrow \hat{k}\geq 1$, with
\begin{equation}\label{eq7}
w(\hat{s}) = \hat{k} q
\end{equation}
and $I(w)\leq\liminf_{n\to\infty}I(w_n) = d$ by weak lower-semicontinuity of\, $I$. In particular, we have $I(w) < \infty$, so that $w\in\Lambda$ and, necessarily,
$\hat{k} > 1$ in \eqref{eq7}. Therefore, $d > 0$ is achieved at $w\in E$.
\endproof
\smallskip

The previous two lemmas show that, for any $S\in \Gamma_2$, we have $\sup_{u\in S}I(u)\geq\inf_{v\in E}I(v)=d$, so that
$$
c_2 \geq d > 0\,.
$$
As the sublevel sets $I^b=\{u\in \Lambda: I(u)\leq b\}$ are complete, standard arguments imply the existence of  a $(PS)$-sequence $(u_n)\subset\Lambda$ at the minimax level $c_2>0$, i.e., there exists $u_n\in\Lambda$ such that
$$
I(u_n) \longrightarrow c_2\ \mbox{ and }\ I^{\prime}(u_n)\longrightarrow 0\,.
$$
In fact, if no such sequence exists, then there exist $\epsilon_0 >0$ and $\delta >0$ such that $||I^{\prime}(u)||\geq \epsilon_0$ if
$|I(u)-c_2|<\delta$. Now it is straightforward to construct a deformation on $\Lambda$ that will deform a set $S\in \Gamma_2$ with $\sup_{u\in S} I(u)<c_2+\epsilon$ (whose existence for $\epsilon$ small follows from the definition of $c_2$) into $\overline{S}\in \Gamma_2$ such that $\sup_{u\in \overline{S}} I(u)<c_2-\epsilon$, contradicting the definition of $c_2$.

Finally, in view of Lemma \ref{lemma1.3}, there exists another $(PS)_{c_2}$-sequence $(v_n)$ and $0\neq v\in\Lambda$ such that $v_n \rightharpoonup v$. Clearly, $v$ satisfies
$$
\left< I^{\prime}(v),\varphi\right> =0\quad\quad \forall\varphi\in C^{\infty}_0 (\Rn)\,,
$$
in other words, $v\in \Lambda\subset H^1$ is a {\it nontrivial} solution of
$$
\ddot{u} + a(t) W^{\prime}(u)=0\,,\quad -\infty < t < \infty\,.
$$
and $\lim_{|t|\rightarrow\infty}v(t)=0$. Furthermore, since $\ddot{v} = - a(t) W^{\prime}(v)$ and $v(t)\longrightarrow 0$ as $t\to\infty$, $(H_2)$ implies that $\ddot{v}\in L^2(\R)$, hence $\dot{v}\in H^1(\R)$. So, an application of Lemma \ref{lemma1.0} (with $v$ replaced by $\dot{v}$) implies $|\dot{v}(t)|\longrightarrow 0$ as $|t|\to\infty$. The proof of Theorem 1.1 is now complete.
\section{Existence of infinitely many homoclinic solutions}          
In this section we take up the question of multiplicity of homoclinic orbits for the singular Hamiltonian system $(SHS)$.

We start by observing that, in view of $(H_2)$, the trivial solution $u=0$ is an isolated critical point of the functional $I$ given in \eqref{functional}. In other words, there exists $r_0 > 0$ such that,
$$
\| u\| \geq r_0\quad\forall\, u\in K\setminus\{ 0\}\,,
$$
where $K$ denotes the set of critical points of $I$. Moreover, one has
\begin{lemma}\label{lemma2.1} For any given $r>0$, it holds that
$$
\inf_{\| u\|\geq r}I(u):=\alpha_r > 0\,.
$$
\end{lemma}
\proof
We argue by contradiction. Assume that there exists a sequence $(u_n)$ such that
$$
I(u_n)\longrightarrow 0\quad\mbox{ with }\quad \| u_n\|\geq r\,.
$$
Since $I$ is coercive by Lemma \ref{lemma1.2}, we have $\| u_n\| \leq R_0$ for some $R_0 > 0$, hence
\begin{equation}\label{eq2.1}
0 < r \leq \| u_n\| \leq R_0\quad\mbox{ and }\quad \| u_n\|_{\infty}\leq\delta\ \forall n\in\mathbb{N}\,,
\end{equation}
for some $\delta > 0$ (recall that $H^1 \subset L^{\infty}$). In particular, we may assume (passing to a subsequence, if necessary) that
$$
\| u_n\|_{\infty}\longrightarrow D \geq 0\,.
$$
We claim that $D > 0$. Otherwise, $(H_2)$ would imply (cf. \eqref{eq1}) $I(u_n)\geq c_0 \| u_n\|^2$ for some $c_0 > 0$ and $n$ large, yielding the contradiction
$$
0 < c_0 r^2 \leq c_0 \| u_n\|^2 \leq I(u_n) \longrightarrow 0\,.
$$
Therefore, we conclude that
$$
\| u_n\|_{\infty}\longrightarrow D > 0\,.
$$
\smallskip

\noindent On the other hand, let $t_n$ be such $|u_n (t_n)| = \| u_n\|_{\infty}$ and note that $I(u_n)\to 0$ implies $\|\dot{u}_n\|_2 \to 0$. Also note that
$$
|u_n (t_n) - u_n (t)| \leq |t_n - t|^{\frac{1}{2}}\|\dot{u}_n\|_2\,,
$$
since $u_n\in H^1$. From this it follows that there exists $T_n \to \infty$ such that
$$
|u_n (t)|\geq \frac{D}{2} > 0\quad\mbox{ for }\ t\in [t_n - T_n, t_n + T_n]\,,
$$
which implies the contradiction
$$
R_0^2 \geq \int_{}|u_n (t)|^2\,dt \geq \int_{[t_n - T_n, t_n + T_n]}|u_n (t)|^2\,dt \longrightarrow \infty\,.
$$
\endproof

\begin{remark}\label{remark2.1}
\begin{description}
\item[(a)]
From Lemma \ref{lemma2.1} and its preceding remark, it follows that there exists $\alpha_0>0$, such that
$$
I(u) \geq \alpha_0\quad\forall\, u\in K\setminus\{ 0\}\,.
$$
\item[(b)] We observe that Lemma \ref{lemma2.1} can be used in place of Lemmas \ref{lemma1.4} and \ref{lemma1.5} to show that $c_2>0$. In fact, if $Cat_{\Lambda}(S)\geq 2$ and $r_1>0$ is chosen such that $||u||_{\infty}<\frac{|q|}{2}$ for all $u\in H^1$ with $||u||\leq r_1$, then $S\cap \left( B_{r_1}(0)\right)^c\not= \emptyset$ (since, otherwise, the homotopy $H$ of Lemma \ref{lemma1.4} would imply $Cat_{\Lambda}(S)=1$) and,  therefore,
$$
c_2=\inf_{S\in \Gamma_2} \sup_{u\in S}I(u)\geq \inf_{||u||\geq r_1}I(u)=\alpha_1>0\,.
$$
\end{description}
\end{remark}
\medskip

Next we consider the structure of $(PS)$-sequences of $I$. Among other things, our proof of infinitely many {\it geometrically distinct} homoclinic solutions for $(SHS)$ (in this case of nonsingular Hamiltonians) will use appropriate versions of a general splitting result for Palais-Smale sequences and of a deformation lemma obtained by Coti-Zelati and Rabinowitz in \cite{9} when considering nonsingular Hamiltonians.

\noindent In order to present these results, let us denote by $\tau_k (v)=v(\cdot - kT)$ the $kT$-shift of $v$, $k\in\mathbb{Z}$. Observe that the periodicity assumption $(A)$ implies that $\tau_k (v)\in K\setminus  \{0\}$ if $v\in K\setminus  \{0\}$. We can make $v$ unique in the class of translations by assuming that $||v||_{\infty}=|v(t_0)|$ for $t_0\in [0,T)$ such that $|v(t)|<|v(t_0)|$ for all $t<t_0$. Such $v$'s are called {\it normalized} critical points. We can now state:

\setcounter{theorem}{0}
\begin{theorem}\label{theorem2.1}(cf. \cite{9}, Proposition 1.24)  Assume $(u_n)$ is a $(PS)_c$-sequence ($c>0$). Then, there exists $l\in\mathbb{N}$, bounded above by a constant depending only on $c$, normalized nonzero critical points $v_1,\ldots,v_l$ of $I$ and corresponding sequences $(k_n^i)\subset\mathbb{Z}$, $1\leq i\leq l$ such that, for a subsequence (still denoted by $(u_n)$), it holds
$$
\| u_n - \sum_{i=1}^l\tau_{k_n^i}v_i \| \longrightarrow 0\,.
$$
\end{theorem}
\proof
Lemma \ref{lemma1.2} implies that any $(PS)_{d}$-sequence is bounded in norm (by a bound depending only on $d$) and, moreover, in view of Lemma \ref{lemma2.1} (see Remark \ref{remark2.1}(a)) we have $I(u)\geq \alpha_0>0$ for all $ u\in K\setminus\{ 0\}$. These two properties are the crucial elements in the proof of Proposition 1.24 of \cite{9} which, as was mentioned above, considers the case of a smooth Hamiltonian. Now it is easily seen that, under assumptions $(A)$ and $(H_1)-(H_4)$, the proof presented in \cite{9} goes through in our case  with obvious modifications.
\endproof

Next, in order to present a version of the deformation lemma that is needed here, we recall a "discreteness" result that was proved in \cite{9} and which is of fundamental importance for the rest of this presentation.

\begin{lemma}\label{finiteness}(\cite{9}, Proposition 1.55)
Let $F\subset \Lambda$ be a finite set of $l\in \N$ points. If
$$
\mathbf{F}_{l}(F):=\{\sum_{i=1}^
{j} \tau_{k_i} (v_i) | 1\leq j\leq l,\ v_i\in F, \ k_i\in \mathbb{Z}\}
$$
then
$$
\mu (F):=\inf\{||x-y|| : x,y\in \mathbf{F}_{l}(F) \}>0\,.
$$
\end{lemma}
We can now state:

\begin{theorem}\label{deformation2.2}(\cite{9}, Proposition 2.2)
Assume, for some $c>0$, that
\begin{itemize}
\item[$(\star)_c$] There exists $\alpha > 0$ such that $I^{c+\alpha}$ has finitely many critical points module $\mathbb{Z}$, i.e., $(K\cap I^{c+\alpha})/ \mathbb{Z}$ is finite\,.
\end{itemize}
Then, for $\bar\epsilon \in (0, \alpha]$ given, there exist $\eta\in C([0,1]\times\Lambda,\Lambda)$, $\epsilon\in (0, \bar\epsilon)$ and $r > 0$, such that
\begin{itemize}
\item[$(i)$] $\eta(0,u)=u$\ $\forall u\in\Lambda$,\ and $\eta(t,\cdot ):\Lambda\rightarrow\Lambda$ is a homeomorphism for any $t\in [0,1]$,
\item[$(ii)$] $I(\eta (s,u))$ is nonincreasing in $s$,
\item[$(iii)$] $\eta (1, I^{c+\epsilon}\setminus N_r(K^{c+\bar\epsilon}_{c-\bar\epsilon})) \subset I^{c-\epsilon}$\
\end{itemize}
where, for $A\subset \Lambda$, we denote $N_r(A)=\{w\in \Lambda : ||w-A||<r\}$. Furthermore $r>0$ can be taken so that
$r<\mu (F)$, with $\mu$ defined as in Lemma \ref{finiteness} for the set $F=(K^{c+\bar\epsilon}_{c-\bar\epsilon})/ \mathbb{Z}:=\{v_1,v_2,\cdot\cdot\cdot v_k\}$.
\end{theorem}

\proof
The proof of  Proposition 2.2 in \cite{9} goes through with no change as  it does not use the explicit form of the functional $I$, relying solely on the splitting of $(PS)_c$-sequences stated in Theorem \ref{theorem2.1}.
\endproof

\medskip

Next, let us recall the definition of the minimax levels $c_m$, $m\in\mathbb{N}$:
$$
c_m := \inf_{S\in\Gamma_m}\sup_{u\in S} I(u)\,,
$$
where
$$
\Gamma_m := \{\, S\subset\Lambda\, |\, S \mbox{ is compact and }Cat_{\Lambda}(S)\geq m\,\}\,.
$$
We have $0=c_1 < c_2 \leq \ldots \leq c_m \leq c_{m+1}\leq \dots<\infty$\ for all $m\in\mathbb{N}$.
The fact that $c_2>0$ was used in section 1 to show existence of one nontrivial homoclinic solution of $(SHS)$.

Finally, we can state and prove our main result:

\begin{theorem}\label{main2.3}
Assume conditions $(A)$, $(H_1)$-$(H_4)$ stated in the Introduction. Then, the singular Hamiltonian system
\vseven
$(SHS)$ \hf
$ \displaystyle \ddot{u} + a(t) W^{\prime}(u)=0\,,\quad -\infty < t < \infty\,,$  \hsf
\vseven
possesses infinitely many {\it geometrically distinct} homoclinic solutions.
\end{theorem}
\proof
We start by noting that, if condition $(\star)_{c_k}$ is satisfied for some $k\in\mathbb{N}$ then $c_k$ is a critical value in view of (the deformation result) Theorem \ref{deformation2.2}. We have two possibilities:
\medskip

\noindent\underline{Case 1}: The levels $c_m$'s are distinct for infinitely many $m$'s.\\
In this case we get infinitely many critical points yielding infinitely many {\it geometrically distinct} homoclinic solutions for $(SHS)$, since either $(\star)_{c_m}$ is satisfied and then $c_m$ is a critical value, or else $(\star)_{c_m}$ is not satisfied, in which case  $(K\cap I^{c_m+\alpha})/ \mathbb{Z}$ is already {\it infinite}.
\medskip

\noindent\underline{Case 2}: There exists $m_0\in\mathbb{N}$ such that $c_{m_0}=c_{m_0 + i}$\ for all $i\geq 1$.\\
In this case we shall assume that $(\star)_{c_{m_0}}$ is satisfied (since, otherwise, there are already infinitely many critical points in $I^{c_{m_0} + \alpha}$, for some $\alpha > 0$, and there is nothing to prove). Therefore, by definition of $c_{m_0 + i}$ and the deformation lemma, there exist $S_i \in \Gamma_{m_0 + i}$ and $\epsilon \in (0, \bar\epsilon)$ (for $\bar\epsilon \in (0,\alpha]$ given) such that
\begin{equation}\label{def_of_c}
\sup_{u\in S_i}\leq c_{m_0} + \epsilon
\end{equation}
and
\begin{equation}\label{eq3i}
\eta (1, S_i) \subset I^{c_{m_0}-\epsilon}\cup \eta (1, N_r(K^{c_{m_0} + \bar\epsilon}_{c_{m_0} - \bar\epsilon}))\,.
\end{equation}
\medskip

Therefore, it follows from \eqref{eq3i} that
$$
\begin{array}{rl}
Cat_{\Lambda}(\eta (1, S_i)) & \leq Cat_{\Lambda}(I^{c_{m_0}-\epsilon}) + Cat_{\Lambda}\eta(1,(N_r(K^{c_{m_0} + \bar\epsilon}_{c_{m_0} - \bar\epsilon})))\\
        & \leq m_0 + Cat_{\Lambda}\eta(1, (N_r(K^{c_{m_0} + \bar\epsilon}_{c_{m_0} - \bar\epsilon})))
\end{array}
$$
On the other hand, since $S_i \in \Gamma_{m_0 + i}$ and $\eta (1,\cdot)$ is a homeomorphism, we also have $\eta (1, S_i) \in \Gamma_{m_0 + i}\,$ and
$$ Cat_{\Lambda}(\eta (1, S_i))= Cat_{\Lambda}( S_i)\geq m_0+i,\quad\quad Cat_{\Lambda}\eta(1, (N_r(K^{c_{m_0} + \bar\epsilon}_{c_{m_0} - \bar\epsilon})))=Cat_{\Lambda}(N_r(K^{c_{m_0} + \bar\epsilon}_{c_{m_0} - \bar\epsilon}))
$$
Hence,
\begin{equation}\label{inf_cat}
Cat_{\Lambda}(N_r(K^{c_{m_0} + \bar\epsilon}_{c_{m_0} - \bar\epsilon})) = \infty\,,
\end{equation}
since $i\in\mathbb{N}$ is arbitrary. In addition, by our assumption, $K^{c_{m_0} + \bar\epsilon}_{c_{m_0} - \bar\epsilon}/ \mathbb{Z}$ is finite, say
$$
K^{c_{m_0} + \bar\epsilon}_{c_{m_0} - \bar\epsilon}/ \mathbb{Z} = \{\, v_1,\ldots, v_l\, \}\,,
$$
where $v_i$'s are normalized critical points. Therefore
$$
K^{c_{m_0} + \bar\epsilon}_{c_{m_0} - \bar\epsilon} = \bigcup_{i=1}^{l}{\cal T}_{\mathbb{Z}}(v_i)\,,
$$
where we are denoting ${\cal T}_{\mathbb{Z}}(v_i) = \{\, v\, |\, v=\tau_p(v_i) \mbox{ for some $p\in\mathbb{Z}$}\,\}$. In view of Lemma \ref{finiteness}, this representation implies that the set
$K^{c_{m_0} + \bar\epsilon}_{c_{m_0} - \bar\epsilon}$ is discrete and, since $r< \mu (\{\, v_1,\ldots, v_l\, \})$, we have
$$
N_r(K^{c_{m_0} + \bar\epsilon}_{c_{m_0} - \bar\epsilon}) = \bigcup_{i=1}^{l}\bigcup_{p\in\mathbb{Z}}B_r (\tau_p (v_i))\,.
$$
\medskip

However, each $\bigcup_{p\in\mathbb{Z}}B_r (\tau_p (v_i))$ has category $1$. Indeed we can first deform the closure of each ball $B_r (\tau_p(v_i))$ to its center $\tau_p(v_i)$ and then use the homotopy
$$
L(s,\tau_p(v_i))(t) = v_i (t + sp)\ ,\quad s\in [0,1]\,,
$$
to deform $\tau_p(v_i)$ to $v_i$ (Clearly this defines a continuous map on $\bigcup_{p\in\mathbb{Z}}B_r (\tau_p (v_i))$ in view of Lemma \ref{finiteness}).
\medskip

 It follows that $Cat_{\Lambda}(N_r(K^{c_{m_0} + \bar\epsilon}_{c_{m_0} - \bar\epsilon}))\leq l$, contradicting \eqref{inf_cat}. Thus, $K^{c_{m_0} + \bar\epsilon}_{c_{m_0} - \bar\epsilon}/ \mathbb{Z}$ is \underline{not finite}, which yields infinitely many {\it geometrically distinct} homoclinic solutions for $(SHS)$. The proof of Theorem \ref{main2.3} is finally complete.
\endproof

\begin{remark}
Note that, when applied to the Lusternik-Schnirelmann level $c_2$, the above argument gives an alternative proof of Theorem \ref{theorem1.1} on existence of one homoclinic solution.
\end{remark}
\section{Appendix}                                                
In this appendix we present a proof of Proposition \ref{Cat}. Although this seems to be a known result, we have not been able to find a reference in the literature. On the other hand, the periodic case is a well known and now classical result. In fact, for $T>0$, let $E=W^{1,2}_{T}(\R,\R^N)$ denote the Sobolev space of $T-$periodic functions from $\R$ to $\R^N$. Let $E_0(q)$ be the subset of $E$ consisting of functions that map into $U:=\R^n\setminus \{q\}$ which, in addition, take a base point
in $[0,T]$ to a base point in $U$. Topologically,  $E_0(q)$ has the same homotopy type as $\Omega(U)$, where for any topological space $X$, $\Omega(X)$ is the space of continuous maps from the unit circle $S^1$ into $X$ which take a base point in $S^1$ to a base point of $X$. It follows from a result of Serre \cite{19} that $Cat_{\Omega(U)}\Omega(U)=\infty$ and, in addition, there are compact subsets of any category $k\in \mathbb{N}$ (see also \cite{10} and \cite{11}). Therefore $Cat_{ E_0(q)}E_0(q)=\infty$ and has the corresponding property as well.
The main idea of our proof is to directly show that $\Lambda$ has the same homotopy type as the set
$$
X=\{ u\in H^1_0(-1,1) : u(t)\not= q, \mbox{ for all } t\in\R\},
$$
which is clearly homotopic equivalent to $E_0(q)$ (with $T=2$). This will in turn show the validity of both parts of Proposition \ref{Cat}.

\proof
Without loss of generality we assume that $|q|=2$.
For $u\in \Lambda$ we define:
$$ T_0(u)=\inf\{ t\in \R^{+} : \int_{|s|\geq t} \left( |u|^2(s)+|\dot u|^2(s)\right) ds\leq \frac{1}{4}\}$$
and the function
$$
\phi_{u}:\R\longrightarrow [0,1],\quad \phi_{u}(t)=\left\{
\begin{array}{ccc}
0 & \mbox{ when } |t| < T_0(u),  \\

|t|-T_0(u) & \mbox{ when }  T_0(u)\leq |t|\leq T_0(u)+1, \\

1 & \mbox{ when }  T_0(u)+1 <|t|.
\end{array}
\right.
$$
In addition, to any $u\in \Lambda$ we associate a positive real number $r(u)$ such that:
$$
v\in H^1, \quad ||v-u||\leq r(u)\Longrightarrow v\in \Lambda,\quad ||v-u||_{\infty}\leq \frac{1}{2}.
$$
The balls $\{ B_{r(u)}(u)\}$ for $u\in \Lambda$ provide an open cover for the (metric) space $\Lambda$ and, therefore, there exists a locally finite refinement $\{W_{\alpha}\}_{\alpha\in I} $ consisting of neighborhoods $W_l\subset B_{r(u_l)}(u_l)$. We choose a continuous partition of unity
$\{P_{\alpha}\}_{\alpha\in I}$ subordinate to   $\{W_{\alpha}\}_{\alpha\in I }$ and finally define
$$
G:\Lambda \longrightarrow X=\{ u\in H^1_0(-1,1) : u(t)\not= q, \mbox{ for all } t\in\R\}\,,
$$
$$
G(u)(t)=u(\frac{t}{1-|t|})\left[ 1-\left( \sum_{\alpha\in I}P_{\alpha}(u)\phi_{u_{\alpha}}( \frac{t}{1-|t|})\right)|t|\right]\,.
$$
Straightforward calculations show that $G$ maps continuously into $H^{1}_0(-1,1)$. The fact that the image of $G(u)$ does not pass through $q$ follows from the following observations:
\begin{description}
\item[1.] If for some $t\in (-1,1)$ (equivalently, for the corresponding $s=\frac{t}{1-|t|}$),   $\sum_{\alpha\in I}P_{\alpha}(u)\phi_{u_{\alpha}}(s)=0$, then $G(u)(t)=u(s)\not= q$ since $u\in\Lambda$.
\item[2.] If for some $t\in (-1,1)$ (equivalently,  $s=\frac{t}{1-|t|}$),   $\sum_{\alpha\in I}P_{\alpha}(u)\phi_{u_{\alpha}}(s)\not=0$, then for some $\alpha \in I$ we have
$$
P_{\alpha}(u)\not= 0,\quad \phi_{u_{\alpha}}(s)\not=0.
$$
Thus $|s|\geq T_0(u_{\alpha})$. But from the definition of $T_0(u_{\alpha})$ and Lemma \ref{lemma1.0} we have $|u_{\alpha}(s)|\leq 1$. And since $P_{\alpha}(u)\not= 0$, we have $ ||u-u_{\alpha}||\leq r(u_{\alpha})$, so that $|u(s)|\leq \frac{3}{2}<|q|$.
Therefore, since $G(u)(t)=u(s).a$ for some $0\leq a\leq 1$, we have $G(u)(t)\not= q$.
\end{description}
Next we consider the map
$$ F: X\longrightarrow \Lambda$$
$$ F(u)(t)=u(\frac{t}{1+|t|})$$
It is easily seen that $F$ is a well defined continuous map. Finally, we calculate the compositions\ $F\circ G$\ and $G\circ F$. We have:
$$
F\circ G: \Lambda\longrightarrow\Lambda, \quad (F\circ G)(u)(t)=u(t)\left[ 1-\left( \sum_{\alpha\in I}P_{\alpha}(u)\phi_{u_{\alpha}}(t)\right)\frac{|t|}{1+|t|}\right]\,,
$$
which is homotopic to the identity map on $\Lambda$ through the homotopy
$$H:[0,1]\times \Lambda\longrightarrow\Lambda,\quad H(\lambda,u)(t)=u(t)\left[ 1-\lambda\left( \sum_{\alpha\in I}P_{\alpha}(u)\phi_{u_{\alpha}}(t)\right)\frac{|t|}{1+|t|}\right]\,.
$$
Similarly, we get
$$
G\circ F: X\longrightarrow X,\quad (G\circ F)(u)(t)=u(t)\left[ 1-\left( \sum_{\alpha\in I}P_{\alpha}(F(u))\phi_{u_{\alpha}}(\frac{t}{1-|t|})\right)|t|\right]\,,
$$
which is homotopic to the identity map on $X$ through the homotopy
$$
\widehat{H}:[0,1]\times X\longrightarrow X,\quad \widehat{H}(\lambda,u)(t)=u(t)\left[ 1-\lambda\left( \sum_{\alpha\in I}P_{\alpha}(F(u))\phi_{u_{\alpha}}(\frac{t}{1-|t|})\right)|t|\right]\,.
$$
Hence $\Lambda$ and $X$ have the same homotopy type and, therefore, the same homotopy invariants. In particular,  $\Lambda$ and $X$ have the same category. This completes the proof of Proposition \ref{Cat}.
\endproof

\noindent \begin{acknowledgement}
We would like to thank Prof.~V.~Benci for a helpful comment regarding Proposition \ref{Cat}.
\end{acknowledgement}
%
%

%

\end{document}